\documentclass[reqno]{amsart} \usepackage{amscd}
\usepackage{epsf}
\usepackage{graphics}
\usepackage{graphicx}
\usepackage{ragged2e}
\usepackage{placeins}
\newtheorem{theorem}{Theorem}[section]

\newtheorem{lemma}[theorem]{Lemma}

\newtheorem{remark}[theorem]{Remark}
\newtheorem{que}[theorem]{Question}
\numberwithin{equation}{section}
\newcommand{\field}[1]{\ensuremath{\mathbb{#1}}}

\newcommand{\EE}{\field{E}}

\usepackage{hyperref}
\hypersetup{
	colorlinks,
	citecolor=red,
	filecolor=black,
	linkcolor=blue,
	urlcolor=black
}

\begin{document}
\title{CMC surfaces of revolution, Elliptic curves, Weierstrass-$\wp$ functions, and Algebraicity}
\author[Rukmini Dey]{Rukmini Dey}
\address{Survey No. 151, Shivakote, Hesaraghatta, Uttarahalli Hobli, Bengaluru, 560089, India.}
\email{rukmini@icts.res.in}
\author[Anantadulal Paul]{Anantadulal Paul}
\address{Survey No. 151, Shivakote, Hesaraghatta, Uttarahalli Hobli, Bengaluru, 560089, India.}
\email{anantadulal.paul@icts.res.in}

\author[Rahul Kumar Singh]{Rahul Kumar Singh}
\address{Indian Institute of Technology Patna,
Bihta, Patna -801106 (Bihar), India.}
\email{rahulks@iitp.ac.in}

\begin{abstract}
This paper establishes an interesting connection between the family of CMC surfaces of revolution in $\mathbb E_1^3$ and some specific families of elliptic curves. As a consequence of this connection, we show in the class of spacelike CMC surfaces of revolution in the $\mathbb E_1^3$, only spacelike cylinders and standard hyperboloids are algebraic. We also show that a similar connection exists between CMC surfaces of revolution in $\mathbb E^3$ and elliptic curves. Further, we use this to reestablish the fact that the only CMC algebraic surfaces of revolution in $\mathbb E^3$ are spheres and right circular cylinders.
\end{abstract}

\maketitle

\section{Introduction}


A surface $S$ in $\mathbb R^3$ is algebraic if it can be expressed as a zero set of some polynomial $P$. Such a surface is given by the set $S:=P^{-1}(0)=\{(x,y,z)\in \mathbb R^3~|~P(x,y,z)=0\}$. $S$ is regular if the gradient of $P$, i.e., $\nabla P(x,y,z)$ is non-zero for all $(x,y,z)\in P^{-1}(0)$. The degree of the polynomial $P$ is called the degree of the surface $S$. It is a known fact that there are many examples of zero mean curvature surfaces (minimal surfaces) in the Euclidean space $\mathbb E^3:=(\mathbb R^3,dx^2+dy^2+dz^2)$ which are algebraic, for instance, the plane, the Enneper surface are algebraic of degree $1$ and $9$ respectively, to name a few. For more such examples, one can check \cite{Nitsche}.

However, in \cite{DoCarmo}, Barbosa and do Carmo showed that in the class of regular non-zero constant mean curvature surfaces (abbreviated as CMC surfaces) in the Euclidean space $\mathbb E^3$, spheres and the right circular cylinders are the only examples of regular algebraic surfaces (see Theorem $1.1$ in \cite{DoCarmo}). A general question in this direction raised by Barbosa and do Carmo is as follows: 
\begin{que}
\label{regular_hypersurface_in_Rn}
What are all regular algebraic hypersurfaces 
with non-zero constant mean curvature in 
$\mathbb{R}^{n+1},~n\geq 2$?
\end{que}
To the best of our knowledge, this question in such a generality is still open. However, these types of 
questions are highly classical and have only recently received much more attention. For instance, 
a few years ago, in \cite{Oscar_Vlad}, the authors prove that there are no CMC algebraic hypersurfaces of degree $3$ in $\mathbb{R}^n, ~n\geq 3$. A partial 
development towards the question \eqref{regular_hypersurface_in_Rn} has already appeared in \cite{Barreto_P}. They show that 
there are no regular algebraic hyper-surfaces with non-zero constant mean curvature in $\mathbb{R}^{n+1}$ for $n \geq 2$, defined 
by odd-degree polynomials. Very recently, in \cite{Globally_sub_analytic_in_R3}, authors have generalised the study of algebraic CMC surfaces into CMC surfaces 
in $\mathbb{R}^3$, which are globally sub-analytic. These subsequent developments play a significant motivation for us to study 
the algebraicity of rotational CMC surfaces. \\

In \cite{DoCarmo}, Barbosa and do Carmo crucially used the structure theorem for CMC surfaces due to Korevaar, Kusner et al. (see \cite{KKS} and \cite{meeks}). Interestingly, they also happen to be the surfaces of revolution in $\mathbb E^3$. One can consider the class of CMC surfaces of revolution in $\mathbb E^3$, which also includes non-regular examples of CMC surfaces like nodoids (they have self-intersection). It was not known if the nodoids are algebraic or not; however, we have shown that they are not algebraic. 
In \cite{DoCarmo}, the authors also show that there are no regular algebraic minimal surfaces in $\mathbb E^3$ except plane (see Proposition $4.1$ in \cite{DoCarmo}). Recall that there are examples of non-regular minimal surfaces which are algebraic. One such example is that of Enneper's surface.  

Based on the above observation, it would be interesting to know all regular algebraic spacelike CMC surfaces (natural analogues of CMC surfaces in $\mathbb E^3$) in the Lorentz-Minkowski space $\mathbb E_1^3:=(\mathbb R^3, dx^2+dy^2-dz^2)$. But no such structure theorem is available for spacelike CMC surfaces in $\mathbb E_1^3$. Hence, we can not adopt the method used in \cite{DoCarmo} to prove the same result in $\mathbb E_1^3$. 

In this paper, one of our aims is to study the algebraicity of spacelike CMC surfaces (may not be regular) of revolution (about a spacelike axis and timelike axis) in $\mathbb E_1^3$ and show that the only algebraic surfaces in this class are spacelike cylinders and hyperboloids. A complete classification of spacelike CMC surfaces of revolution in $\mathbb E_1^3$ exists (see Theorem $4$ in \cite{IF}).
There are in total four main steps (in order), which we have used here to prove our claim:
\begin{itemize}
\item Following \cite{kenmotsu}, we can obtain  Kenmotsu's type parametrization for spacelike CMC surfaces of revolution in $\mathbb E_1^3$ (this is also obtained in \cite{IF}).	
\item For a given surface of revolution, one can associate an elliptic curve using the above parametrization.
\item For every elliptic curve whose discriminant is non-zero, there exists a Weierstrass -$\wp$ function, a well-known fact (see  \cite{LW}).
	\item A Weierstrass-$\wp$ function does not satisfy any algebraic relation.
\end{itemize}

We also study the algebraicity of CMC surfaces of revolution in $\mathbb E^3$. By following exactly the same steps as above, we prove that unduloids and nodoids (non-regular CMC surfaces of revolution in $\mathbb E^3$) are not algebraic. We also reproduce the result in \cite{DoCarmo} that only regular algebraic CMC surfaces of revolution are right circular cylinders and spheres. It is important to note that in \cite{DoCarmo}, the authors do not consider the question of algebraicity of non-regular CMC surfaces.

Our paper is organised as follows: In section \ref{SCMCR}, we obtain the formula for parametrizations for a family of generating curves of CMC surfaces of revolution in $\mathbb E_{1}^3$. 
In section \ref{UNW}, we show that the generating curves for CMC surfaces of revolution in $\EE_1^{3}$ give rise to a family of elliptic curves (see Lemma \ref{CMC_ELL} and Lemma \ref{CMC_ELL_1}). Further, by using the fact that for every non-singular elliptic curve (discriminant non-zero) there is an associated Weierstrass- $\wp$ function and then by applying the method of contradiction, we conclude that spacelike cylinders and hyperboloids are the only examples of spacelike CMC surfaces of revolution (about spacelike axis and timelike axis) in $\mathbb E_1^3$ which are algebraic (see Theorem \ref{MT} and Theorem \ref{MT1}). Finally, in Section \ref{UNEE}, using the same idea as in Section \ref {UNW}, we establish the fact that in the class of non-zero CMC surfaces of revolution in $\mathbb E^3$, spheres and right circular cylinders are the only examples which are algebraic (see Lemma \ref{CMC_ELL_2} and Theorem \ref{MT2}).

\section{Spacelike CMC surfaces of revolution in $\mathbb E_1^3$}\label{SCMCR}
We denote a smooth arc length parametrised planar spacelike curve in $xz$-plane of $\mathbb E_1^3$ by $\alpha(s)=(x(s),~0,~ z(s))$, where $s\in I$, $I$ an open interval, and $x(s) > 0$ for all $s\in I$. Then, the surface of revolution generated by the curve $\alpha(s)$ about $x$-axis (spacelike) has the following parametric representation
\begin{equation}\label{SX}
S(s, \theta) = (x(s), z(s)\sinh \theta, z(s)\cosh \theta),\end{equation} $s \in I$, $\theta \in \mathbb R$,
and the parametrization of the surface of revolution generated by the curve $\alpha(s)$ about $z$-axis (timelike) is given by
\begin{equation}\label{TX}
	T(s, \theta) = (x(s)\cos \theta, x(s)\sin \theta, z(s)),
	\end{equation} $s \in I$, $\theta \in (0, 2\pi)$. 

All our surfaces of revolution here are assumed to have constant mean curvature $H$; such surfaces are classified by a non-negative real parameter $B$ (see the upcoming subsections).
\subsection{Spacelike CMC surfaces of revolution about a spacelike axis}\label{SLS}
We start with the parametrization \eqref{SX}. Then, its first and second fundamental forms are given by
 $  I_{X}=ds^2+z(s)^2d\theta^2$ and $II_{X}=(x'(s)z''(s)-z'(s)x''(s))ds^2+x'(s)z(s)d\theta^2$ respectively. Now, using these expressions of fundamental forms, we compute its mean curvature $H$ as follows:
 $$H=\frac{-x'(s)z(s)-z(s)^2(x'(s)z''(s)-z'(s)x''(s))}{2z(s)^2}.$$
Equivalently, $H$ satisfy the following ODE:
\begin{equation}\label{HODE}
    2Hz(s)+x'(s)+z(s)(x'(s)z''(s)-z'(s)x''(s))=0.
\end{equation}
Multiplying \eqref{HODE} by $x'(s)$, then using $x'(s)^2-z'(s)^2=1$ and $x'(s)x''(s)-z'(s)z''(s)=0$, \eqref{HODE} reduce to the following ODE:
\begin{equation}\label{HODE1}
    2Hz(s)x'(s)+(z(s)z'(s))'+1=0.
\end{equation}
Similarly, after multiplying \eqref{HODE} by $z'(s)$, we get
\begin{equation}\label{HODE2}
    2Hz(s)z'(s)+(x'(s)z(s))'=0.
\end{equation}
Before proceeding further, let us recall a split-complex number denoted by $t=x+ky$, where $x,~y$ are two real numbers. The conjugate of $t$ is denoted by $t^{\ast} = x - k y$. Since $k^2 = 1$, the absolute value of $t$ is given by $\mid t \mid = \sqrt{x^2-y^2}$. Also, the analogue of Euler's formula in split-complex number system has the form
\[ e^{k\theta} = \cosh\theta + k \sinh\theta.  \]
Now, if we set $Y(s)=z(s)z'(s)+kz(s)x'(s)$. Then, by combining \eqref{HODE1} and \eqref{HODE2} we get a first-order linear (split-complex) ODE:
\begin{equation}\label{MoDE}
    Y'(s)+2kHY(s)+1=0.
\end{equation}
When $H=0$, the equation \eqref{MoDE} reduces to 
\begin{equation}\label{MoDE1}
    Y'(s)+1=0.
\end{equation}
Solving \eqref{MoDE1}, we check that the parametrization $(x(s),z(s))$ of the generating curve satisfies the following equation $z=c\cos \frac{x}{c}$, and when rotated, we get a maximal surface, and it is obvious that it is not algebraic. \\

When $H\neq 0$, the general solution of \eqref{MoDE} is given by 
\begin{equation}\label{MODES}
    Y(s)=\frac{1}{2kH}((1+2kHC)-e^{2kHs})e^{-2kHs}.
\end{equation}
Let $Be^{k\theta}=1+2kHC$, where $\theta \in \mathbb R$, $B\in \mathbb R$, and $C$ a split complex number. Then \eqref{MODES} changes to 
\begin{equation}\label{MODES1}
    Y(s)=\frac{Be^{k(\theta-2Hs)}-1}{2kH}.
\end{equation}

Assuming $H>0$ and by doing a translation of the arc-length parameter $s$, we can ignore $\theta$ in \eqref{MODES1}. Hence, we can solve for $x(s)$ and $z(s)$ easily, obtaining
\begin{multline}\label{UNS}
 S(s;H,B):=(x(s), z(s)) \\ = \left(\int_0^s \frac{B \cosh 2Ht - 1}{\sqrt{1+B^2-2B \cosh 2Ht}}dt, \frac{1}{2H} \sqrt{1+B^2-2B \cosh 2Hs}\right).
\end{multline}
When rotated about $x$-axis (a spacelike axis) gives a spacelike CMC surface of revolution.
\subsection{Spacelike CMC surfaces of revolution about a timelike axis}
By following exactly a similar computation, as done in the previous subsection \ref{SLS}, gives us the following parametrization of the generating curve lying in  $xz$-plane,  given by 
\begin{multline}\label{UNT}
T(s;H,B):= (x(s),z(s)) \\ = \left( \frac{1}{2H} \sqrt{2B \cosh 2Hs -B^2-1}, ~~ \int_0^s \frac{1-B \cosh 2Ht}{\sqrt{2B \cosh 2Ht-B^2-1}}dt\right).
\end{multline}
When rotated about $z$-axis (a timelike axis) gives a spacelike CMC surface of revolution.
For the sake of completeness, we give an overview of the computation; we begin with the parametrization \eqref{TX} of a spacelike surface of revolution about $z$-axis (timelike). Then, its first and second fundamental forms are given by
$  I_{X}=ds^2+x(s)^2d\theta^2$ and $II_{X}=(x'(s)z''(s)-z'(s)x''(s))ds^2+x(s)z'(s)d\theta^2$ respectively. We also obtain the following ODE satisfied by $H$ 
\begin{equation}\label{HODET}
	2Hx(s)+z'(s)+x(s)(x'(s)z''(s)-z'(s)x''(s))=0.
\end{equation}
Again, by multiplying the ODE \eqref{HODET} by $z'(s)$ and $x'(s)$ and using the fact that $x'(s)^2-z'(s)^2=1$, we obtain the following two simultaneous ODEs, given by

\begin{equation}\label{HODET1}
	2Hx(s)z'(s)+(x(s)x'(s))'-1=0.
\end{equation}
and 
\begin{equation}\label{HODET2}
	2Hx(s)x'(s)+(z'(s)x(s))'=0.
\end{equation}
Now, by setting $W(s)=x(s)x'(s)+kx(s)z'(s)$, we see that the ODEs \eqref{HODET1} and \eqref{HODET2} is equivalent to the following first order linear (split-complex) ODE:
\begin{equation}\label{ODEW}
	W'(s)+2kHW(s)-1=0.
\end{equation}

When $H=0$, the solution of the ODE \eqref{ODEW} gives us the elliptic catenary $x=c\sinh \frac{z}{c}$, and when rotated about $z$-axis gives us the elliptic catenoid $x^2+y^2=c^2\sinh^{2}\frac{z}{c}$, which is a maximal surface (see \cite{kobayashi}) and it is obvious that it is not algebraic. When $H\neq 0$, the solutions of \eqref{ODEW} are given by \eqref{UNT}.
\subsection{Unduloids and Nodoids}\label{Undnod}
The curves $S(s;H,B)$ and $T(s;H,B)$ given by the equations \eqref{UNS} and \eqref{UNT},  give parametrizations for spacelike \textit{unduloid} ($0<B<1$) and spacelike \textit{nodoid} ($B>1$), when rotated about $x$-axis and $z$-axis respectively. They are the natural analogues of \textit{nodoids} and \textit{unduloids} in $\mathbb {E}^3$ (see \cite{kenmotsu}). 
 When $B=0$, it corresponds to a spacelike cylinder; when $B=1$, it corresponds to a standard hyperboloid. 

The spacelike curves $S(s;H,0)$ and $S(s;H,1)$ represent the line $z=\frac{1}{2H}$ and the circle $x^2-z^2=\frac{-1}{H^2}$ in $\mathbb L^3$, when rotated about $x$-axis, generates the spacelike cylinder (given by $z^2-y^2=\frac{1}{4H^2}$) and the standard hyperboloid ($x^2+y^2-z^2=\frac{-1}{H^2};\; z>0$), with non-zero constant mean curvature $H$ respectively and it is also clear that they are algebraic. 

Again, we check that for $B=0$ and $B=1$, i.e., the spacelike curves $T(s;H,0)$ and $T(s;H,1)$, generates an imaginary spacelike cylinder (given by $x^2+y^2=\frac{-1}{4H^2}$) and the standard hyperboloid ($x^2+y^2-z^2=\frac{-1}{H^2};\;x>0$), with non-zero constant mean curvature $H$ respectively, when rotated about $z$-axis and it is obvious that they are algebraic.

\section{Spacelike CMC surfaces of revolution in $\mathbb E_1^3$ and Elliptic curves}\label{UNW}

\subsection{Spacelike unduloids, nodoids having timelike axis of revolution and associated Weierstrass-$\wp$ functions}

We begin with the following lemma
\begin{lemma}\label{CMC_ELL}
    The generating curves (namely undulary and nodary, see \eqref{UNT}) for a family of spacelike non-zero CMC surfaces of revolution in $\mathbb E_1^{3}$ (about $z$-axis) identify a family of elliptic curves.
\end{lemma}
\begin{proof}
From  $\eqref{UNT}$, we have $x(s) = \frac{1}{2H} \sqrt{2B \cosh 2Hs -B^2-1}$ and
\begin{align}
\label{z-in-terms-of-B}
z(s)=&\int_0^s \frac{1-B \cosh 2Ht}{\sqrt{2B \cosh 2Ht-B^2-1}}dt.
\end{align}  
Substituting $u = \cosh 2Ht$ we get
\begin{equation}
\label{z-transfm-in-u}
z(s) = \frac{1}{2H} \int_1^{\cosh 2Hs} \frac{1- Bu}{\sqrt{(B^2+1) - 2Bu - (B^2+1)u^2 + 2 Bu^3}}.
\end{equation}
Now set $w = u+c$, we have
\begin{align}
\{(B^2+1) - 2Bu - (B^2+1)u^2 + 2 Bu^3\}^{1/2}&=\{(B^2+1) + 2Bc - (B^2+1)c^{2} - 2Bc^{3} \nonumber \\
& + (6Bc^2 + 2c(B^2+1) - 2B)w \nonumber \\
& - ((B^2+1) + 6Bc)w^2 + 2Bw^3\}^{1/2}.
\label{eq:Lr1}
\end{align}
Next we choose $c = - \frac{B^2+1}{6B}$ then 
the right hand side of Eq. (\ref{z-transfm-in-u}) will reduce to 
\begin{equation}
\sqrt{(B^2+1) + 2Bc - (B^2+1)c^{2} - 2Bc^{3} + (6Bc^2 + 2c(B^2+1) - 2B)w  + 2Bw^3}.
\end{equation}
Putting 
\begin{align}
l& = (B^2+1) + 2Bc - (B^2+1)c^{2} - 2Bc^{3}; \nonumber \\
m& = 6Bc^2 + 2c(B^2+1) - 2B; \nonumber \\
n& = 2B;
\end{align}
then equation (\ref{z-transfm-in-u}) will reduce to
\begin{equation}
z(s)=\frac{1}{2H}\int_{1+c}^{\cosh 2Hs + c}\frac{(1+Bc) - Bw}{\sqrt{l + mw + nw^3}}dw.
\label{z-transfm-in-w}
\end{equation}
and letting
\begin{equation}
v^2=l+mw+nw^3
\label{eq:Ly2}
\end{equation}
followed by the transformation $ w\rightarrow \sqrt[3]{\frac{4}{n}} \tilde{w} $, equation (\ref{eq:Ly2}) will reduce to
\begin{equation}\label{EC}
v^2=4 \tilde{w}^3-(-m \sqrt[3]{\frac{4}{n}})\tilde{w} -(-l).
\end{equation}
The equation \eqref{EC} represents a family of elliptic curves parametrized by $B$, where $B\in (0,1)$ corresponds to an \textit{undulary}  and $B\in (1, \infty)$ corresponds to \textit{nodary}.

\end{proof}

Now we prove that the spacelike CMC surfaces \textit{unduloids} and \textit{nodoids} (about a timelike axis) are not algebraic. 
\begin{theorem}
\label{MT}
    The spacelike CMC surfaces of revolution in $\EE_{1}^3$ (whose axis of revolution is a timelike axis) for  
  $B\in (0, \infty)\setminus \{1 \}$  cannot be algebraic except spacelike cylinders and standard hyperboloids.
\end{theorem}
\begin{proof}
  By Lemma \ref{CMC_ELL}, we know that the discriminant of an elliptic curve \eqref{EC} is given by $(-m \sqrt[3]{\frac{4}{n}})^3-27(-l)^2 $, when evaluated this gives us a degree $8$ polynomial $2(B^2 - 1)^4$
  and we note that
\begin{equation}
(-m \sqrt[3]{\frac{4}{n}})^3-27(-l)^2 \neq 0; ~\text{except for}~ B =  1.
\end{equation}


Hence, we conclude that the discriminant of the elliptic curve \eqref{EC} can not be zero for any  $B\in(0, \infty) \setminus \{1 \}$.
Therefore, there is a lattice $L$ such that
\begin{equation}
g_2(L)=-m \sqrt[3]{\frac{4}{n}};\ \ g_3(L)=-l
\end{equation} 
\noindent
and corresponding to this, there will be a Weierstrass-$\wp$ function, which will satisfy
\begin{equation}\label{WP}
(\wp '(\tilde{z}))^2=4(\wp(\tilde{z}))^3-g_2\wp(\tilde{z})-g_3,
\end{equation}
\noindent
where $\tilde{z}\in \mathbb C$ (for details see ~\cite[page 269]{LW}).
Note that we can always choose a neighbourhood of $\tilde{z}$, sufficiently small, such that in this neighbourhood $\wp'(\tilde{z})\neq 0$, and hence we can differentiate the above equation to obtain
\begin{equation}
2\wp''(\tilde{z})=12\wp^2(\tilde{z})-g_2.
\end{equation}
Now, if we replace $w$ by $\wp(\tilde{z})$ in equation (\ref{z-transfm-in-w}) and using the fact that the Weierstrass-$\wp$ function takes the real values only on the real and imaginary axes, the new limits of integration will be changed to $\tilde{z}_0=\wp^{-1}(1+c)$ and $\tilde{z} = \wp^{-1}( \cosh 2 H s + c) $ ($\tilde{z}$ depends on $s$), where we can choose $\tilde{z}_0$ to be any one of the two values of $\wp^{-1}(1+c)$, we choose the real $\tilde{z}_0$. Next, we select a neighbourhood of $\tilde{z}_0$ small enough such that $\wp(\tilde{z})$ restricted to the portion of real axis contained in this neighbourhood, becomes real-valued. Since we have chosen $\tilde{z}_0$ to be real, $\tilde{z}$ can be taken to be real, and hence we obtain a real-valued integral as $\wp(\tilde{z})$ is real. 
Now, if we set $p=1+Bc$; $q=-B$ and then using equation \eqref{WP} in (\ref{z-transfm-in-w}) (note that here $v$ will be replaced by $\wp'(z)$),  we get
\begin{align}
z(t)&=\frac{1}{2H}\int_{t_0}^{t}(q\wp(t)+p)dt \nonumber \\
&=\frac{p}{2H}(t-t_0) + \frac{q}{2H}\int_{t_0}^{t}\wp(t)dt,
\end{align}
where we have changed $\tilde{z}_0$ to $t_0$ and $\tilde{z}$ to $t$.
We also have
\begin{align}\label{XT}
x(t)&=  \frac{1}{2H} \sqrt{2B \cosh 2Hs -B^2-1}\nonumber \\ 
&= \frac{1}{2H} \sqrt{2B(\wp(t)-c)-B^2-1}. \\
\textnormal{Therefore},~x(t)^2&= \frac{-(B^2+2Bc+1) + 2B\wp(t) }{4 H^2}\nonumber \\
&= {c_1+c_2 \wp(t)},
\label{eq:Lr2}
\end{align}
where $c_1 = \frac{-(B^2+2Bc+1)}{4 H^2}$ and $c_2 = \frac{2B}{4 H^2}$.

Therefore, in this new parameter $t$, we can describe the generating curve as $T(t, H, B) = \left(x(t), z(t) \right)$ and the corresponding surface of revolution it generates is expressed as $(x_1, x_2, x_3) = (x(t) \cos(\theta),~ x(t) \sin(\theta),~ z(t))$. We know that any spacelike surface of revolution about $x_3$-axis (timelike axis) is given by zeroes of an equation of the form $x_1^2+x_2^2-r(x_3)$ (see \cite{Sneddon-book}). Our claim is that the function $r(x_3)$ cannot be a polynomial. On the contrary, assume that $r(x_3)$ is a polynomial.

Since we have $r(x_3)=x_1^2+x_2^2 = x(t)^2$; $x(t) = \sqrt{r(x_3)}$. Therefore, from \eqref{XT}, we get
\begin{equation}\label{RX3}
r(x_3) = c_1 + c_2\wp(t).
\end{equation}
Set $\alpha=\frac{p}{2H}$, $\beta=\frac{q}{2H}$, $ c_1 = \frac{-(B^2+2Bc+1)}{4 H^2}$, and
$c_2=\frac{2B}{4H^2}$.
Note that all constants appearing here ultimately depend on $B$. Also, note that $x_3=z(t)$ and $t=f(x_3)$, where $f$ is the inverse of $z$. This imply $\frac{dx_3}{dt}\frac{dt}{dx_3}=1$ and $\frac{dt}{dx_3}=\frac{1}{\alpha+\beta\wp(t)}$.\\

Differentiating equation \eqref{RX3} w.r.t $x_3$ we get
\begin{align}
\frac{dr}{dx_3}&=c_2\wp'(t)\frac{dt}{dx_3} \nonumber \\
&=c_2\wp'(t)\frac{1}{\alpha+\beta\wp(t)} \nonumber \\
&=R_1(\wp)\wp'.
\label{eq:r6}
\end{align}
 
where $R_1(\wp)$ is a rational function of $\wp$.
Again
\begin{align}
\frac{d^2r}{dx_3^2}&=c_2\frac{d}{dx_3}\left( \frac{\wp'}{\alpha+\beta\wp}\right) \nonumber \\
&=c_2\left\lbrace \frac{\wp''}{\alpha+\beta\wp}-\beta\left(\frac{\wp'}{\alpha+\beta\wp}\right)^2 \right\rbrace \frac{dt}{dx_3}
\label{eq:r7}
\end{align}
Now, using the relations
$(\wp ')^2=4\wp^3-g_2\wp-g_3$ and $\wp''=6\wp^2-\frac{g_2}{2}$ in Eq. (\ref{eq:r7}) we get
\begin{align}
\frac{d^2r}{dx_3^2}&=c_2\left\lbrace\frac{12\wp^2-g_2}{2(\alpha+\beta\wp)^2}-\beta\frac{4\wp^3-g_2\wp-g_3}{(\alpha+\beta\wp)^3}\right\rbrace \nonumber \\
&=R_2(\wp),
\end{align}
where $R_2(\wp)$ is again a rational function of $\wp$; continuing in this way and by frequently using the above relations, the derivatives at the odd stage and even stage will have the form
\begin{equation}
\frac{d^{2k-1}r}{dx_3^{2k-1}}=R_{2k-1}(\wp)\wp',
\end{equation}
\begin{equation}
\frac{d^{2k}r}{dx_3^{2k}}=R_{2k}(\wp)
\end{equation}
respectively. Observe that on the left-hand side of the above equations, we are repeatedly taking the derivative of a polynomial $r(x_3)$ (by assumption); hence, after a while, there exists a natural number (say $m$) such that
$\frac{d^{m}r}{dx_3^{m}}=0$, which can be either even or odd depending on the degree of $r(x_3)$ ($m$ will be odd if the degree is even and vice-versa). 
If $m$ is even, i.e., $m=2k$
\begin{equation}
\frac{d^{2k}r}{dx_3^{2k}}=0 \Longrightarrow R_{2k}(\wp)=0.
\end{equation}
This tells us that $\wp$ satisfies an algebraic relation, which is not possible as we know $\wp$ is transcendental.
Next if $m$ is odd, i.e., $m=2k-1$
\begin{equation}
\frac{d^{2k-1}r}{dx_1^{2k-1}}=0 \Longrightarrow R_{2k-1}(\wp)\wp'=0
\label{eq:r8}
\end{equation}
Now observe that we can always select a neighbourhood around the point $z_0$ so there will be no zeroes of $\wp'$ in that neighbourhood. Then equation (\ref{eq:r8}) will imply $R_{2k-1}(\wp)=0$. So again, we are getting an algebraic relation satisfied by $\wp$. Thus, in any case, a polynomial equation is satisfied by $\wp$, which is impossible.
Thus, we arrive at a contradiction. Hence, we have shown that the spacelike CMC surfaces of revolution (unduloids and nodoids) about a timelike axis cannot be algebraic.
\end{proof}

\textbf{Note:} 
When the discriminant of the corresponding elliptic curve vanishes, namely when $B = 1$, we cannot get the corresponding Weierstrass-$\wp$ function. However, when $B = 1$, we have already shown that the corresponding surface is part of the standard hyperboloid, which is anyway algebraic (see the subsection \ref{Undnod}).  

\subsection{Spacelike unduloids, nodoids having a spacelike axis of revolution and associated Weierstrass-$\wp$ functions}
Following exactly the same steps as in the proof of Lemma \ref{CMC_ELL} and Theorem \ref{MT}, we show that the CMC surfaces of revolution (about a spacelike axis), are not algebraic for $B\in (0,\infty)\setminus \{1\}$. For $B=1$, the corresponding surface of revolution is part of the standard hyperboloid, which is algebraic (see the subsection \ref{Undnod}).
\begin{lemma}\label{CMC_ELL_1}
    The generating curves (namely undulary and nodary, see \eqref{UNS})  for a family of spacelike non-zero CMC surfaces of revolution in $\mathbb E_1^{3}$ (about $x$-axis) identify a family of elliptic curves.
\end{lemma}

\begin{proof}
Continuing as in the proof of Lemma \ref{CMC_ELL}, we will derive a cubic equation parametrized by a non-zero constant $B$. Let us start with the parametrization
\[ x(s) = \int_0^s \frac{B \cosh 2Ht - 1}{\sqrt{1+B^2-2B \cosh 2Ht}}dt. \]
Next, proceeding as in Lemma \ref{CMC_ELL}, we substitute $\cosh 2Ht = \tilde{u}$, $\tilde{w} = \tilde{u} + \tilde{c}$, and by choosing $\tilde{c} = -\frac{1+B^2}{6B}$, we arrive at
\begin{equation}
x(s)=\frac{1}{2H}\int_{1+\tilde{c}}^{\cosh 2H s + \tilde{c}}\frac{B\tilde{w}-(1+B\tilde{c})}{\sqrt{\tilde{l} + \tilde{m} \tilde{w} + \tilde{n} \tilde{w}^3}}d\tilde{w},
\label{z-transfm-in-w-xAxis}
\end{equation}
where

\begin{align}
\tilde{l}&= -(B^2+1) -2B \tilde{c} + (B^2+1) \tilde{c}^{2} + 2B\tilde{c}^{3}; \nonumber \\
\tilde{m}&= 2B - 2 \tilde{c}(B^2 + 1) - 6B \tilde{c}^{2}; \nonumber \\
\tilde{n}&= -2B.
\end{align}
Hence, the lemma follows as before.

\end{proof}

\begin{theorem}
\label{MT1}
    The spacelike CMC surfaces of revolution in $\EE_{1}^3$ (whose axis of revolution is a spacelike axis) for $B\in (0,\infty)\setminus \{1\}$, cannot be algebraic except spacelike cylinders and standard hyperboloids.
\end{theorem}

\section{Spacelike CMC surfaces of revolution in $\mathbb E^3$ and Elliptic curves}\label{UNEE}
By following the ideas developed in section \ref{UNW}, in this section, we show that the only algebraic CMC rotational surfaces of revolution in $\mathbb E^3$ are the sphere and the right circular cylinders. This also aligns with the fact that the only algebraic non-zero CMC regular surfaces are spheres and cylinders (see \cite{DoCarmo}).  

 As earlier, we consider the curve $\alpha(s)=(x(s),~ y(s))$, where $s\in I$, $I$ an open interval, and $y(s) > 0$ for all $s\in I$. We consider the surface of revolution generated by the curve $\alpha(s)$, which has a local parametric representation
$$S(s, \theta) = (x(s), y(s) cos(\theta), y(s) sin(\theta) ) = (x_1, x_2, x_3),$$ $s \in I$, $\theta \in (0, 2\pi)$, and $x_1$,$x_2$, $x_3$ denotes Cartesian coordinates of $\mathbb R^3$.
Here, we consider the surfaces of revolution as having mean curvature $H$ constant.
Such surfaces are classified by a non-negative real parameter $B$ (see \cite{kenmotsu}). For instance, $B=0$ corresponds to a sphere, a surface of revolution with a mean curvature constant, and $B=1$ corresponds to a cylinder, which are already known to be algebraic. Since each $B\in (0,\infty)$ corresponds to a surface of revolution having constant mean curvature, we denote such a surface by $S(H,~B)$, whose generating curve will be denoted by $X(s; H, B)$. We let $X(s;H,B) = (x(s), y(s))$. Then by Proposition $3.2$ in \cite{kenmotsu},  we have
\begin{equation}
X(s;H,B)=\left(\int_0^s \frac{1+B \sin2Ht}{\sqrt{1+B^2+2B \sin2Ht}}dt, \frac{1}{2H} \sqrt{1+B^2+2B \sin2Hs}\right),
\end{equation}
where $B \in [0,\infty)$.\\ 

Similar to the Lorentzian case, we can prove the following lemma.

\begin{lemma}\label{CMC_ELL_2}
    The generating curves for a family of nonzero CMC surfaces of revolution in $\mathbb E^3$ identify a family of elliptic curves.
\end{lemma}

\begin{remark}
    Interestingly, the discriminant of the associated elliptic curves does not vanish for any $B\in (0, \infty)$. 
\end{remark}
By following exactly the same idea as in Section \ref{UNW}, we can show that the only CMC surfaces of revolution, which are also algebraic, are spheres and cylinders. 
\begin{theorem}\label{MT2}
    The nonzero CMC surfaces of revolution in $\mathbb E^3$, corresponding to $0<B<1$ (unduloids) and $B > 1$ (nodoids), are not algebraic.
\end{theorem}


\section{Acknowledgment}

Rukmini Dey would like to acknowledge the support of Department of Atomic Energy, Government of India, under project no. RTI4001.  Rahul Kumar Singh would like to acknowledge the external grant he has obtained, namely MATRICS (File No. MTR/2023/000990) that has been sanctioned by the SERB.
\bibliography{CMC}
\bibliographystyle{siam}
\end{document}